\documentclass[10pt]{article}
\usepackage{amssymb,latexsym,amsmath,epsfig,amsthm} 

\makeatletter

\renewcommand\section{\@startsection {section}{1}{\z@}
{-30pt \@plus -1ex \@minus -.2ex}
{2.3ex \@plus.2ex}
{\normalfont\normalsize\bfseries}}

\renewcommand\subsection{\@startsection{subsection}{2}{\z@}
{-3.25ex\@plus -1ex \@minus -.2ex}
{1.5ex \@plus .2ex}
{\normalfont\normalsize\bfseries}}

\renewcommand{\@seccntformat}[1]{\csname the#1\endcsname. }

\makeatother

\newtheorem{theorem}{Theorem}
\newtheorem{lemma}{Lemma}
\newtheorem*{conjecture*}{Conjecture}

\newtheorem*{corollary*}{Corollary}
\theoremstyle{remark}
\newtheorem*{remark*}{Remark}
\newtheorem*{definition}{Definition}
\begin{document}
\begin{center}
\uppercase{\bf Two divisibility problems on subset sums}
\vskip 20pt
{\bf Konstantinos Gaitanas}\\
{\tt kostasgaitanas@gmail.com}\\
\vskip 10pt
\end{center}
\vskip 30pt
\vskip 30pt

\centerline{\bf Abstract}
\noindent
We consider two problems regarding some divisibility properties of the subset sums of a set $A\subseteq \{1, 2, \ldots ,n\}$. At the beginning, we study the cardinality of $A$ which has the following property: For every $d\le n$ there is a non empty set $A_d\subseteq A$ such that the sum of the elements of $A_d$ is a multiple of $d$. Next, we turn our attention to another problem: If all subset sums of $A$ form a multiple free-sequence, what can we say about the structure of $A$? We give some asymptotics for the first problem and improve some already existing results for the second one.
\pagestyle{myheadings}  
\thispagestyle{empty} 
\baselineskip=12.875pt 
\vskip 30pt
\section{Introduction} 
Since the beginning of the 20-th century, there has been a great interest in problems related to sums representable in the form $\sum\limits_{a\in A}\epsilon_a\cdot a$ where $\epsilon_a\in \{0, 1\}$. These are called the \textit{subset sums} of $A$. In particular, there are two kind of questions frequently asked. The first one is how the subset sums are distributed in $\mathbb{Z}_n$  and the second one is, if $A\subseteq \{1, 2, \ldots ,n\}$ how the subset sums of $A$ are distributed in $\{1, 2, \ldots ,n \}$. We refer the reader to \cite{3} for further discussion of the literature and additional references. Regarding the second case, a well- known open problem exists:
 If $A$ has distinct subset sums (SSD), what is the maximum cardinality of $A$? Erd\H{o}s conjectured that $|A|=\frac{\ln n}{\ln 2}+O\left(1\right)$ and Guy and Conway showed \cite{2} that $|A|\ge\lfloor \log_2 {n}\rfloor+2$ where $\log_2 n$ denotes the base 2 logarithm of $n$. At this point,  Erd\H{o}s's conjecture remains open. However, in \cite{1} a new question is posed which is analogous to the SSD conjecture. In this analogy, ``distinct subset sums'' corresponds to ``subsets sums not dividing one another '', where the authors prove that $\frac{\log n}{\log 2}-1<|A|<\frac{\log n}{\log 2}+\frac{\log \log n}{2\log 2}+c$. Surprisingly, it seems that this question did not draw much attention since then.  We deal with this problem at the second part of the paper. Firstly, we introduce a problem which seems that it has not been investigated so far: If for every $d\le n$ there is $\sum\limits_{a\in A}\epsilon_a\cdot a$ such that $d\mid \sum\limits_{a\in A}\epsilon_a\cdot a$, how small can $|A|$ be? For the latter two cases we obtain some partial results. 
\section{Preliminaries}
Throughout this paper, we shall use the standard abbreviation $[n]$, for the set of natural numbers from $1$ to $n$ and $S(X)$ for the sum of the elements of a set $X$. In addition, we will need the following definition.
\begin{definition}
Let $A\subseteq [n]$. We say that $A$ is a multiple of $[n]$ , if for every $d\in N$ there is a nonempty $A_d\subseteq A$ such that $d\mid S(A_d)$.
\end{definition}
\section{Problem one}
It is not very hard to construct an explicit set $A$ which is a multiple of $[n]$ for every $n$, since every natural number can be represented as a sum of powers of $2$ (and hence, divide itself). For example, we see that $A=\{1,2,4,\ldots, 2^z\}$ is a multiple of $[2^z]$, since it contains all the powers of $2$ not greater than $2^z$. What we can see however is that even for small values of $n$ it is possible that $|A|<\lfloor \log_2{n}\rfloor+1$. The natural question is then to ask how small can $|A|$ be with respect to $n$. In order to obtain a non-trivial bound, we will need the following inequality interesting on its own:\\
\begin{lemma}
For every $n\ge2$ the inequality $n>\sum\limits_{d\mid n} \frac{d}{\tau(d)}$ holds.
\begin{proof}
Here $\tau(d)$ denotes the number of positive divisors of $d$. It is easy to see that the function $\frac{n}{\tau(n)}$ is multiplicative. If $n=\prod\limits_{i=1}^{r}p_i^{a_i}$ is the prime factorization of $n$, the sum in the lemma can be written as $\prod\limits_{i=1}^{r}\sum\limits_{j=0}^{a_j}\frac{p_j^{a_j}}{\tau(p_j^{a_j})}$, which since $\tau(p_j^{a_j})=a_j+1$ is equal to $\prod\limits_{i=1}^{r}\sum\limits_{j=0}^{a_j}\frac{p_j^{a_j}}{a_j+1}$. But $\sum\limits_{j=0}^{a_j}\frac{p_j^{a_j}}{a_j+1}<1+p_j^1+\ldots+p_j^{a_j-1}+\frac{p_j^{a_j}}{a_j+1}=\frac{p_j^{a_j}-1}{p_j-1}+\frac{p_j^{a_j}}{a_j+1}$ which is less than $p_j^{a_j}$ since $a_j\ge 1$. We conclude that $\prod\limits_{i=1}^{r}\sum\limits_{j=0}^{a_j}\frac{p_j^{a_j}}{a_j+1}<\prod\limits_{i=1}^{r}p_i^{a_i}=n$.
\end{proof}
\end{lemma}
This lemma will be usefull in order to prove the following theorem. We make use of the $O$ notation in the standard way. We say that $f(x)=O\left(g(x)\right)$ if $g(x)>0$ and there is a $c>0$ such that $|f(x)|\le cg(x)$ for sufficiently large $x$.
\begin{theorem}
Let $[n]=\{1,2,\ldots ,n\}$ and  $l(n)$ be the smallest possible cardinality of a set $A$ which is a multiple of $[n]$. Then $l(n)= \log_2{n}-O\left( \log_2\log {n}\right)$.
\begin{proof}
It is a routine matter to prove that every natural number not greater than $n$ divides some even number not greater than $2n$, which are represented by the subset sums of $A=\{2^1, 2^2,\ldots ,2^z\}$ where $2^z\le n<2^{z+1}$. This shows that $A$ is a multiple of $[n]$ and hence $l(n)\le \log_2{n}$.\\
We can get a lower bound for $l(n)$ by using the fact that if $d\mid x$ then $d\le x$. Suppose that $A$ is a multiple of $[n]$ and let $\sum\limits_{i=1}^{2^{l(n)}-1}S(A_i)$ be the sum of all non-empty subset sums of $A$. A simple counting argument shows that this is equal to $\sum\limits_{i=1}^{2^{l(n)}-1}S(A_i)=2^{l(n)-1}(a_1+\ldots+a_{l(n)})<2^{l(n)-1}l(n)\cdot n$, since every element of $A$ is at most $n$. We restrict attention on subset sums which are divisible by numbers $D$,  such that $\frac{n}{2}< D\le n$ since for every $d\in [n]$ corresponds a number $D$ such that $d\mid D$. If $D\mid S(A_i)$ then $S(A_i)\ge D$ and then (using the previous lemma), $S(A_i)>\sum\limits_{d\mid D} \frac{d}{\tau(d)}$. This shows that $\sum\limits_{i=1}^{2^{l(n)}-1}S(A_i)>\sum\limits_{d\le n} \frac{d}{\tau(d)}$.\\
It is proved in \cite{4} that $\sum\limits_{n\le x} \frac{1}{\tau(n)}=\frac{x}{\sqrt{\log x}}\left (c_0+\frac{c_1}{\log x}+\ldots +\frac{c_v}{(\log x)^v}+O\left( \frac{1}{(\log x)^{v+1}}\right)\right)$ where the $c_0, c_1,\ldots ,c_v$ are constants. Using partial summation we can find that $\sum\limits_{n\le x} \frac{n}{\tau(n)}>\frac{c_0x^2}{2\sqrt{\log x}}$ if $x$ is sufficiently large.\\
From this we get the inequality $2^{l(n)-1}l(n)\cdot n>\frac{c_0n^2}{2\sqrt{\log n}}$ and if we take logarithms with base $2$ we can find that $l(n)>\log_2n-\log_2l(n)-\frac{1}{2}\log_2(\log n)+\log_2c_0$.\\
Since $l(n)\le \log_2n$, the right hand side is equal to $\log_2n-O\left(\log_2(\log n)\right)$. This completes the proof.\\
\end{proof}
\end{theorem}
One could ask whether it is possible to sharpen the conclusion of Theorem 1 further, to assert that $l(n)= \log_2{n}-O\left(1\right)$. However this question is likely to be almost as hard to settle as the distinct subset sum conjecture, and out of reach of the methods of this paper.
\section{Problem two}
We now turn our attention to the other problem. Let $A=\{a_1, \ldots ,a_k\}\subseteq [n]$. It is stated in \cite{1} as property R that the subset sums of $A$ form a multiple-free sequence (no one divides another). The authors discuss what is the maximum value of $|A|$ such that property R holds. In particular, if $n=2^z, z\in \mathbb{N}$ they show using elementary methods that property R holds for $A=\{2^z-1, 2^z-2^1, 2^z-2^2, \ldots ,2^z-2^{z-2}, 2^z-2^{z-1}\}$ and thus $|A|\ge z$. It seems that since then this bound has not been improved. However, numerical data suggests that the following seems to be true:
\begin{conjecture*}
Let $A=\{a_1, \ldots ,a_{k}\}\subseteq [2^z]$. If for every $i, j\le2^k-1$ the quotient $\frac{S(A_i)}{S(A_j)}$ is not a power of $2$, then $k=z+O\left(1\right)$.
\end{conjecture*}
If someone is more optimistic, he can expect that actually if $k=z+1$, the quotient of the two subset sums is exactly $2$. We say that $A=\{a_1, \ldots ,a_k\}\subseteq [2^z]$ possesses property R* if for every $i, j\le2^{k}-1$ the quotient $\frac{S(A_i)}{S(A_j)}$ is not a power of $2$. Property R*  immediately shows the subset sums must be different. Indeed, if $S(A_i)=S(A_j)$, then $\frac{S(A_i)}{S(A_j)}=1=2^0$. In order to support the optimistic side of the conjecture, we observe that equal subset sums create a quotient which is exactly $2$. Indeed, if $S(A_i)=S(A_j)$, then assuming that $A_i, A_j$ are disjoint (otherwise the same summands cancel out) we get that $S(A_i)+S(A_j)=2S(A_j)$. But $A_i\cup A_j=A_m$ for some $m\neq i, j$ and so, $S(A_m)=2S(A_j)$. Our goal from now on is to carry further results in order to understand better the structure of these type of sets. We begin with the following observation:
\begin{theorem}
Let $A=\{2^z-1, 2^z-2^1, 2^z-2^2, \ldots ,2^z-2^{z-2}, 2^z-2^{z-1}\}$. For every $c\le 2^z$, property R (and thus, R*) does not hold true for $A\cup\{c\}$.
\begin{proof}
To see this, let $c=2^{b_1}+\ldots +2^{b_x}$ be it's binary representation, where $0\le b_x<\ldots <b_1$. We observe that $c+(2^z-2^{b_1+1})+\ldots +(2^z-2^{b_x+1})=(2^z-2^{b_1})+\ldots +(2^z-2^{b_x})$ where we may assume that both sides of the equation do not share a same summand. Note that every term inside the parentheses is an element of $A$ and thus of $A\cup\{c\}$. This means that two subset sums of $A\cup\{c\}$ are equal which is absurd. 
\end{proof}
\end{theorem} 
In order to prove the following theorem, we will need to make use of a celebrated result in a version not so well-known. For the convenience of the reader we also give a short proof. \\
\begin{lemma}
If we choose a sequence of $v=\lfloor\frac{x+1}{2}\rfloor+1$ integers not greater than $x$, then there are two elements of the sequence such that their quotient is a power of $2$.
\begin{proof}
Every element of the sequence can be written in the form $a_i=2^{b_i}\cdot w_i, i\le v$ where $w_i$ is odd.There are at most $v-1$ odd numbers not greater than $x$, so using the pidgeonhole principle, there must be $i, j$ such that $w_i=w_j$ and $b_i\le b_j$ which shows that $\frac{a_j}{a_i}=2^{b_j-b_i}$.
\end{proof}
\end{lemma}
Using this result we prove the following:
\begin{theorem}
If $A$ possesses property R*, then $a_1+\ldots +a_k\ge 2^{k+1}-3$.
\begin{proof}
We have $2^k-1$ non-empty subset sums from $A$ which are all not greater than $a_1+\ldots+a_k$. For the shake of contradiction, suppose that $a_1+\ldots+a_k\le 2^{k+1}-4$. Taking $x=2^{k+1}-4$, we get $2^k-1=\lfloor\frac{x+1}{2}\rfloor+1$ and using the previous lemma we arrive at a contradiction.
\end{proof}
\end{theorem}
We continue by showing that a set which possesses property R* must satisfy an additional condition. By \textit{subset differences} we denote all sums of the form $\sum\limits_{a\in A}\epsilon_a\cdot a$ where $\epsilon_a\in \{-1, 0, 1\}$ and at least one of the $\epsilon_a$'s is equal to $-1$. It is not difficult to show that the number of positive (but not necessarily different) subset differences of $A$ is equal to $\frac{3^k+1}{2}-2^k$. To see this, we take all possible sums of the form $\sum\limits_{a\in A}\delta_a\cdot a$ where $\delta_a\in \{-1, 0, 1\}$. We exclude the empty sum and observe that for every positive sum of this form, it's additive inverse (which is negative) is also of this form.We also wish to exclude all the non-empty subset sums which are $2^k-1$ at number. This shows that the number of positive subsset differences is $\frac{3^k-1}{2}-(2^k-1)=\frac{3^k+1}{2}-2^k$. We view the next result as further evidence towards the truth of our conjecture.
\begin{theorem}
If $A=\{a_1, \ldots,a_k\}$  possesses property R*, then no subset difference is equal to any subset sum.
\begin{proof}
Suppose $S(A_i)=S(A_j)-S(A_l)$, where we may assume that $A_i\cap A_j=\emptyset$ and $A_j\cap A_l=\emptyset$. We rewrite the equation as $S(A_i)+S(A_l)=S(A_j)$. We observe that if $A_i\cap A_l=\emptyset$, then we have two subset sums which are equal and as we proved earlier, this is a contradiction. If $A_i\cap A_l\neq \emptyset$, we get $S(A_i\setminus A_l)+S(A_l\setminus A_i)+2S(A_i\cap A_l)=S(A_j)$. We add $S(A_i\setminus A_l)+S(A_l\setminus A_i)$ on both sides of the equation and since all sets are disjoint, we see that this is equivalent to $2S(A_p)=S(A_q)$, where $S(A_p)=S(A_i\setminus A_l)+S(A_l\setminus A_i)+S(A_i\cap A_l)$ and $S(A_q)=S(A_i\setminus A_l)+S(A_l\setminus A_i)+S(A_j)$. This completes the proof.
\end{proof}
\end{theorem}
The above results of course do not prove the mentioned conjecture. Nevertheless, these demonstrate that the set $A$ must have special structure which possibly cannot exist if $|A|=k$ is much bigger than $z$.
\makeatletter
\renewcommand{\@biblabel}[1]{[#1]\hfill} 
\makeatother

\end{document}